\documentclass[a4paper,10pt]{amsart}

\usepackage[centertags]{amsmath}
\usepackage{amsfonts}
\usepackage{amssymb}
\usepackage{amsthm}
\usepackage{newlfont}
\usepackage{euscript}

\newtheorem{thm}{Theorem}
\newtheorem{lem}{Lemma}
\newtheorem{prop}{Proposition}
\theoremstyle{definition}
\newtheorem{dfn}{Definition}
\newtheorem{ex}{Example}
\newtheorem{rem}{Remark}
\numberwithin{equation}{section}

\newcommand{\brakes}[1]{\left(#1\right)}

\newcommand{\eqdef}{\stackrel{\mathrm{def}}{=}}

\newcommand{\lims}{\underset{n\rightarrow\infty}{{\lim\sup}}\,}
\newcommand{\norm}[1]{\left\Vert#1\right\Vert}
\newcommand{\abs}[1]{\left\vert#1\right\vert}
\newcommand{\figbrakes}[1]{\left\{#1\right\}}
\newcommand{\sqrbrakes}[1]{\left[#1\right]}

\newcommand{\sign}{\mathrm{sign}\,}
\renewcommand{\Re}{{\Bbb R}}
\newcommand{\CC}{{\Bbb C}}
\newcommand{\DD}{{\Bbb D}}
\newcommand{\NN}{{\Bbb N}}
\newcommand{\ZZ}{{\Bbb Z}}
\newcommand{\TT}{{\Bbb T}}
\newcommand{\Xf}{{\frak X}}

\newcommand{\Gc}{\mathcal{G}}
\newcommand{\Fc}{\mathcal{F}}
\newcommand{\Nc}{\mathcal{N}}
\newcommand{\Lc}{\mathcal{L}}
\newcommand{\1}{1\!\!{\mathrm I}}
\newcommand{\ax}{\Re^+}

\newcommand{\eps}{\varepsilon}
\newcommand{\be}{\begin{equation}}
\newcommand{\ee}{\end{equation}}
\newcommand{\eqd}{\mathop{=}\limits^{d}}

\begin{document}
\title{Invariance principle for additive functionals of Markov chains}
\author{Yuri N.Kartashov, Alexey M.Kulik}
\abstract{We consider a sequence of additive functionals
$\{\phi_n\}$, set on a sequence of Markov chains $\{X_n\}$ that
weakly converges to a Markov process $X$. We give sufficient
condition for such a sequence to converge in distribution,
formulated in terms of  the characteristics of the additive
functionals, and  related to the Dynkin's theorem on the
convergence of $W$-functionals.  As an application of the main
theorem, the general sufficient condition for convergence of
additive functionals in terms of transition probabilities of the
chains $X_n$ is proved.}\endabstract



\email{kulik@imath.kiev.ua} \maketitle

\section{Introduction}

Let a sequence of processes $X_n=X_n(\cdot)$ be given, converging
in distribution (in some sense, e.g., in a sense of convergence of
finite-dimensional distributions, distributions in spaces $\CC$ or
$\DD$, etc.) to a limit process $X=X(\cdot)$. Also let the family
of functionals $\phi_n$ of the processes $X_n$ be given. Assume
that they are additive in an appropriate sense  with respect to
time variable. The general question, considered in the present
paper, is  what an information about the limit behavior of the
distributions of functionals $\phi_n$ can be obtained in a
situation where the processes $X_n, X$ possess certain Markov
properties. The starting
 point in our considerations is provided by the comparatively simple, but important
particular case of the problem outlined above, in which all the
processes $X_n$ coincide. In this situation, $\phi_n$ are a
functionals of the same process $X$, and if $X$ is Markov process
and $\phi_n$ are $W$-functionals (see \cite{dynkin}, Chapter 6),
then  their limit behavior, according to the well known theorem by
E.B.Dynkin (\cite{dynkin}, Theorem 6.4), is determined by the
limit behavior of their characteristics (that is, their
 expectations).

In the present paper we consider the processes $X_n$ that differ
one from another. The class of sequences of processes $X_n$,
considered in the framework of our approach, contains sequences of
Markov chains with appropriately normalized time, embedded into
$\CC$ or $\DD$ (for example, by means of standard operations of
linearization or construction of graduated processes), and weakly
convergent to Markov process $X$. Important partial case  is
provided by random broken lines (or random step functions) $X_n$,
constructed by a random walk in $\Re^d$ and weakly convergent to a
homogenous stable process $X$ (particularly, to the Brownian
motion).

 We show that, under some structural assumption about processes
 $X_n,X$ (the condition is that  the sequence
$\{X_n\}$ provides \emph{Markov approximation} for the process
$X$), the full analogue of the Dynkin's theorem takes place: if
the characteristics of functionals $\phi_n$ converge weakly to the
characteristics of $W$-functional $\phi$ of the limit process $X$,
then the distributions of $\phi_n$ converge to the  distribution
of $\phi$. Our method of proof is based on $L_2$-estimates for the
distance between additive functionals, similar to those given in
Lemma 6.5 \cite{dynkin}. The proof of these estimates is concerned
with a preliminary construction of processes $X_n,X$ on  one
probability space in such a way, that the functionals
$\phi_n,\phi$, associated initially with a  different processes,
are interpreted as a functionals of one two-component process. The
(some kind of) Markov property of the two-component process is
essential for the estimates, analogous to those given in Lemma 6.5
\cite{dynkin}; the structural assumption mentioned above is just
the claim for such a property to hold true in an appropriate form.

The method, proposed by authors, allows one to reduce the problem
of studying of asymptotic behavior of the distributions of
additive functionals to \emph{a priori} more simple problem of
studying of their means. In our opinion, it provides a good
addition to the available methods of studying the limit behavior
of additive functionals both for the important partial case of
random walks (we do not give  the detailed review here, referring
the reader to monographs
\cite{skor_slob},\cite{bor_ibrag},\cite{revez}, papers
\cite{bass_khoshn92},\cite{Shi_Che_yor} and reviews there), and
for general Markov chains. Among the latter, it is necessary to
mention the method that is based on  the passing to the limit in
the difference equations that describe characteristic functions of
additive functionals of Markov chains, and ascends to the works of
I.I.Gikhman at 50-ies (see \cite{gikh_1},\cite{gikh_2}, also
\cite{portenko_68} and the survey paper \cite{portenko_98}).

The structure of the article is following. In Chapter 2, we
introduce the notion of Markov approximation and give examples
that illustrate it. In Chapter 3, the main theorem of the article
is introduced and proved. In Chapters 4,5,  the two elementary
examples of application of this theorem  are given. In Chapter  6,
the main theorem is applied to the proof of a general sufficient
condition for weak convergence of additive functionals, set on the
sequence of Markov chains, that is formulated in terms of
transition probabilities of the chains.

\section{Markov approximation.}

Further we assume that the processes $X_n, X$ are defined on $\ax$
and have a locally compact metric phase space $(\Xf,\rho)$. We say
that the  process $X$ possesses the Markov property at the time
moment $s\in\ax$ w.r.t. filtration $\{\mathcal{G}_t,t\in\ax\}$, if
  $X$ is adapted to this filtration and for each $k\in \NN,
t_1,\dots,t_k>s$ there exists a stochastic kernel $\{P_{st_1\dots
t_k}(x,A), x\in \Xf, A\in\mathcal{B}(\Xf^k)\}$ such
that\be\label{20}
E[\1_{A}((X(t_1),\dots,X(t_k)))|\Gc_s]=P_{st_1\dots t_k}(X(s),A)
\quad \hbox{a.s.,} \quad A\in\mathcal{B}(\Xf^k). \ee The measure
$P_{st_1\dots t_k}(x,\cdot)$ has a natural interpretation as the
finite-dimensional distribution of $X$ at the points
$t_1,\dots,t_k$, conditioned by $\{X(s)=x\}$;  we denote below
$P_{st_1\dots
t_k}(x,\cdot)=P((X(t_1),\dots,X(t_k))\in\cdot|X(s)=x)$.

\begin{rem}\label{r21} In some cases, (\ref{20}) implies
the following functional analogue of (\ref{20}): \be\label{201}
E[\1_{\cdot}(X|_s^\infty)|\Gc_s]=E[\1_{\cdot}(X|_s^\infty)|X(s)],
\ee where $X|_s^\infty$ denotes the trajectory of the process $X$
on the time interval $[s,+\infty)$, considered as an element of
appropriate functional space. For instance, if the Kolmogorov's
sufficient condition for existence of continuous modification
holds true both for unconditional and conditional distributions of
$X$, then (\ref{201}) holds with $X|_s^\infty$ considered as an
element of $C([s,+\infty),\Xf)$.
\end{rem}

Everywhere below we assume that the process $X$ possesses the
Markov property w.r.t. its canonic filtration at every point $s\in
\ax$ and for the processes $X_n$ the same property holds true at
every point of the type ${i\over n}, i\in \ZZ_{+}$ (the choice of
the denominator here is quite arbitrary; it is possible to put
  any expression $N(n)\to \infty, n\to \infty$ instead of $n$, but
we avoid to do this in order to  shorten the notation).

The next definition is introduced in \cite{kulik}.

\begin{dfn}\label{d11} The sequence $\{X_n\}$
provides Markov approximation for the process $X$, if for
arbitrary $\gamma>0, T<+\infty$ there exists $K(\gamma,T)\in \NN$
and a sequence of two-componential processes $\{\hat Y_n=(\hat
X_n,\hat X^n)\}$, defined on another probability space, such that

(i) $\hat X_n\mathop{=}\limits^d X_n, \hat X^n \mathop{=}\limits^d
X$;

(ii) the process $\hat Y_n$, together with the processes $\hat
X_n, \hat X^n$, possesses the Markov property at the points ${i
K(\gamma,T)\over n}, i\in \NN$ w.r.t. filtration $\{\hat
{\Fc}_t^n=\sigma(\hat Y_n(s),s\leq t)\};$

(iii)
$\mathop{\lim\sup}\limits_{n\to+\infty}P\left(\sup\limits_{i\leq
{Tn\over K(\gamma,T)}} \rho\left(\hat X_n\left({i K(\gamma,T)\over
n}\right),\hat X^n\left({i K(\gamma,T)\over
n}\right)\right)>\gamma\right)<\gamma.$
\end{dfn}

\begin{rem} Condition (ii) implies that,
 for $i,k\in \NN, t_1,\dots,t_k>{i K(\gamma,T)\over n}, (x,y)\in
\Xf^2$, the marginal distributions $P\brakes{(\hat
Y_n(t_1),\dots,\hat Y_n (t_k))\in\cdot|\hat Y_{n}({i
K(\gamma,T)\over n})=(x,y)}$ are equal to
$P\brakes{(X_n(t_1),\dots,X_n(t_k))\in\cdot |X_n({i
K(\gamma,T)\over n})=x}$ and
$P\brakes{(X(t_1),\dots,X(t_k))\in\cdot |X\brakes{i
K(\gamma,T)\over n}=y}$ respectively.
\end{rem}

Let us give some examples that illustrate Definition \ref{d11}.

\begin{ex} Let $\{\xi_k\}$ be a sequence of i.i.d random vectors
in $\Re^d$ with $E\|\xi_k\|^{2+\delta}_{\Re^d}<+\infty$ for some
$\delta>0$. Assume $\{\xi_k\}$ to have zero mean and identity for
their covariance matrix. Let us introduce the sequence of
processes $X_n$ ("random broken lines") on $\ax$ by  \be\label{21}
X_n(t)={S_{k-1}\over \sqrt{n}}+(nt-{k+1})\left[{S_{k}\over
\sqrt{n}}-{S_{k-1}\over \sqrt{n}}\right],\quad t\in\left[{k-1\over
n},{k\over n}\right),\quad k\in \NN,\ee where
$S_n=\sum_{k=1}^n\xi_k$. Then $X_n$ converge by distribution in
$C(\ax,\Re^d)$ to the Brownian motion $X$ in $\Re^d$.

It is shown in \cite{kulik} that the sequence $\{X_n\}$ provides
Markov approximation for the process $X$ (part I. of Theorem 1
\cite{kulik}). On the other hand, in the same paper (part II. of
the same Theorem) the following effect is revealed. Let us denote
by $\mathbf{K}(\gamma,T)$ the minimal constant $K(\gamma,T)$ such
that there exists a process $\hat Y_n$ satisfying conditions
(i)-(iii) of Definition \ref{d11}. Then, in all the cases except
one trivial case $\xi_k\sim \Nc(0,I)$, for each fixed $T>0$  the
convergence $\mathbf{K}(\gamma,T)\to +\infty, \gamma\to 0+$ takes
place. In other words, while the accuracy of approximation of the
Brownian motion $X$ by the random walk $X_n$ becomes better (this
accuracy is described by the parameter $\gamma$), the Markov
properties of the pair of processes $(X,X_n)$ necessarily become
worse (these properties are characterized by
$\mathbf{K}(\gamma,T)$).
\end{ex}

\begin{ex} Let $\{\xi_k\}$ be i.i.d random variables, belonging to
the normal domain of attraction for $\alpha$-stable distribution
$\Lc$, $\alpha\in (0,2)$. By the definition, this means that
$$
n^{-{1\over\alpha}}[S_n-a_n]\Rightarrow {\Lc},\quad
a_n=\begin{cases}
0,&\alpha\in(0,1)\\
nE \xi_1,&\alpha\in(1,2)\\
n^2E\sin {\xi_1\over n},&\alpha=1\end{cases}
$$
(\cite{feller}, Chapter XVII.5). In order to shorten the notation,
we assume that $a_n\equiv 0$ and consider processes $X_n$ on $\ax$
of the type \be\label{22}
X_n(t)=n^{-{1\over\alpha}}{S_{k-1}}+(nt-{k+1})\left[n^{-{1\over\alpha}}{S_{k}}-
n^{-{1\over\alpha}}{S_{k-1}}\right],\quad t\in\left[{k-1\over
n},{k\over n}\right),\quad k\in \NN.\ee Then $X_n$ converge by
distribution in $\DD(\ax)$ to the homogeneous process with
independent increments $X$ in $\Re$, for which $X(1)-X(0)\eqd \Lc$
(we call such process a process an $\alpha$-stable one).

It is shown in \cite{kulik} (Theorem 2) that the sequence
$\{X_n\}$ provides Markov approximation for the process $X$.
Furthermore, in this situation, on the contrary to the previous
example, $\mathbf{K}(\gamma,T)=1$ for all $\gamma,T$. This means
that, in this case, the Markov properties do not become worse
while accuracy of approximation improves.
\end{ex}

\begin{rem}\label{r13} The last example shows that the property
of Markov approximation does imply, in general, the convergence of
distributions of the processes $X_n$ to the distribution of $X$ in
$\CC=C(\ax,\Xf)$ even if $X_n$ has continuous trajectories. The
same can be said about convergence in $\DD=\DD(\ax,\Xf)$ (we
 omit the corresponding example).
\end{rem}

 Let us remark that the approach, introduced in the present
paper, is closely related to the
 Skorokhod's method of embedding of random walk into Wiener process by means of
 of appropriate sequence of stopping moments (\cite{skor_issledovania}), widely
 used in literature.
 The basic idea is the same: we have to construct
 two processes on the same probability space, with the pair
 keeping Markov or martingale properties. However, the Skorokhod's method,
 while being quite efficient for one-dimensional random walks that
 approximate Wiener process, is much less appropriate in a
multi-dimensional situation or for stable domain of attraction.
Examples 1 and 2 show  that the claim for the Markov approximation
to hold true is not restrictive, at least for all basic classes of
random walks with no regard to the dimension of the phase space or
to the type of limit distribution.

The following example shows that the property of Markov
approximation is "stable" in the following sense. This property is
preserved under construction  of a new pair $(Z_n,Z)$ from the
pair $(X_n,X)$, possessing this property, in some regular way
(e.g., as a solution of a family of stochastic equations).

\begin{ex}\label{e3} Let $X_n,X$ be as in Example 1,
functions $a:\Re^m\to \Re^m, b:\Re^d\to \Re^{d\times m}$ be
Lipschitz and $b^*(x)b(x)>0, x\in \Re^m$ (the sign $^*$ denotes
the operation of taking of the adjoint matrix). Define
\be\label{23} Z_n\brakes{{k+1\over n}}=Z_n\brakes{{k\over
n}}+a\brakes{Z_n\brakes{{k\over n}}}{1\over
n}+b\brakes{Z_n\brakes{{k\over n}}}\Delta X_n\brakes{{k\over
n}},\quad Z_n(0)=z, \ee
 $\Delta X_n({k\over n})\equiv [X_n({k+1\over
n})-X_n({k\over n})]$. Then (\cite{jacod_shir},
\cite{kurtz_protter}) $Z_n$ converge by distribution in
$C(\ax,\Re^m)$ to the process $Z$, defined by SDE \be\label{24}
dZ(t)=a(Z(t))dt+b(Z(t))dX(t),\quad Z(0)=z, \ee
 where $X$ is the Brownian motion in $\Re^d$. It is natural to call the sequence $Z_n$
the difference approximation of the diffusion process $Z$.

Let us show that the sequence $\{Z_n\}$ provides Markov
approximation for the process $Z$. For arbitrary $\gamma,T$, we
construct a pair $(\hat X_n, \hat X^n)$, corresponding to
processes $X_n, X$ and satisfying conditions of Definition
\ref{d11} (such construction is possible due to  Example 1).

Let us construct the processes $\hat Z_n, \hat Z^n$ as the
functionals of the processes $\hat X_n, \hat X^n$ by equalities
(\ref{23}),(\ref{24}) with $X_n$ replaced by $\hat X_n$ and $X$
replaced by $\hat X^n$ (note that (\ref{24}) has unique strong
solution, hence this  procedure is correct). By the construction,
the pair $(\hat Z_n, \hat Z^n)$ satisfies condition (i) of
Definition \ref{d11}. It is easy to verify that the Markov
condition (ii) for the pair $(\hat X_n, \hat X^n)$ holds in the
functional form (\ref{201}) with $\hat Y_n|_s^\infty$ considered
as an element of $C([s,+\infty),\Re^d\times \Re^d)$ (see Remark
\ref{r21}). Hence, the pair $(\hat Z_n, \hat Z^n)$ also satisfies
condition (ii) of Definition \ref{d11}. Let us write
$$\Delta(\gamma)=\mathop{\lim\sup}\limits_{n\to+\infty}P\left(\sup\limits_{i\leq
{Tn\over K(\gamma,T)}} \rho\left(\hat Z_n\left({i K(\gamma,T)\over
n}\right),\hat Z^n\left({i K(\gamma,T)\over
n}\right)\right)>\gamma\right), $$ and show that
\be\label{25}\Delta(\gamma)\to 0+,\quad \gamma\to 0+.\ee Note that
(\ref{25}) immediately implies Markov approximation: for arbitrary
$\delta>0$ we chose, using (\ref{25}), $\gamma=\gamma(\delta)$
such that inequalities  $\gamma<\delta$ and
$\Delta(\gamma)<\delta$ hold. Then the pair $(\hat Z_n, \hat
Z^n)$, constructed by the scheme described above, satisfy
Definition \ref{d11} with the constant $\gamma$ replaced by
$\delta$ (note that, under this construction, the value
$K(\delta,T)\equiv K_Z(\delta,T)$ for the pair $(\hat Z_n, \hat
Z^n)$ is expressed through the same value for the pair $(\hat X_n,
\hat X^n)$ by $K_Z(\delta,T)=K_X(\gamma(\delta),T)$).

Now assume that (\ref{25}) does not hold, then there exist
constant $c>0$ and sequence $\gamma_k\to 0+, n_k\to +\infty$ such
that \be\label{26} {K(\gamma_k,T)\over n_k}\to 0, \quad
P\left(\sup\limits_{i\leq {Tn_k\over K(\gamma_k,T)}}
\rho\left(\hat Z_n\left({i K(\gamma_k,T)\over n_k}\right),\hat
Z^n\left({i K(\gamma_k,T)\over
n_k}\right)\right)>\gamma_k\right)>c.\ee Consider the sequence of
four-component processes $(\hat X_{n_k},\hat X^{n_k},\hat
Z_{n_k},\hat Z^{n_k})$. Every component of this sequence is weakly
compact in $C(\ax,\Re^d)$ or $C(\ax,\Re^m)$, hence the whole
sequence is also weakly compact in $C(\ax,\Re^d\times\Re^d\times
\Re^m\times\Re^m)$. Consider arbitrary limit point $(\hat
X_{*},\hat X^{*},\hat Z_{*},\hat Z^{*})$ (in a sense of
convergence by distribution) of this sequence. It follows from
(\ref{26}) that \be\label{27} P(Z_*\not=Z^*)>0.\ee It follows from
Theorem 2.2 \cite{kurtz_protter} (see also Chapter 9.5
\cite{jacod_shir}) that the processes $Z_*,Z^*$ satisfy SDE
(\ref{24}) with $X$ replaced by $X_*,X^*$. However, the SDE
(\ref{24}) possesses the property of pathwise uniqueness (see
\cite{yamada_vatanabe}), and the property (iii) of the  pair
$(X_{n_k},\hat X^{n_k})$ implies that the processes $X_*, X^*$
coincide a.s. Therefore, the processes $Z_*, Z^*$ also coincide
a.s., that contradicts to (\ref{27}) and show that our assumption
that $\Delta(\gamma)\not\to 0+, \gamma\to 0+$ is false.
\end{ex}

The examples given above show that the claim for the  Markov
approximation to hold is not very restrictive, and is provided in
a typical situations. On the other hand, this claim is strong
enough to provide one the opportunity to  obtain an analog of the
Dynkin's theorem; this will be shown in the next chapter.

\section{Main theorem}
We consider the functionals of the  type
\be\label{31}\phi_{n}^{s,t}(Y)\eqdef \sum_{k:s\leq
k/n<t}F_{n}\brakes{Y\brakes{\frac{k}{n}},Y\brakes{\frac{k+1}{n}},
\dots,Y\brakes{\frac{k+L-1}{n}}},\quad 0\leq s<t,\ee where the
functions $F_{n}(\cdot)$ are nonnegative, $L$ is a fixed integer.
Together with the functionals $\phi_n$, that are "stepwise"
functions  w.r.t.  every time variable, we consider random broken
lines, related to these functions:
$$
\psi_n^{s,t}=\phi_{n}^{{j-1\over n}, {k-1\over
n}}+(ns-j+1)\phi_{n}^{{j-1\over n}, {j\over
n}}+(nt-k+1)\phi_{n}^{{k-1\over n}, {k\over n}},\quad
s\in\left[{j-1\over n},{j\over n}\right),t\in\left[{k-1\over
n},{k\over n}\right).
$$
We interpret the random broken lines $\psi_n$ as a random elements
in space $C(\TT,\Re^+)$, where $\TT\eqdef\{(s,t)|0\leq s\leq t\}$.

If process $Y$ possesses Markov property w.r.t. the filtration,
associated with this process, at the points of the type $s={i\over
n}, i\in \ZZ_+$, then, for functional $\phi_n$, its characteristic
$f_n$ is naturally defined by the formula \be\label{32}
f_n^{s,t}(x) \eqdef E[\phi_{n}^{s,t}(Y)|Y(s)=x],\quad s={i\over
n}, i\in \ZZ_+, t>s, x\in \Xf. \ee Note that the functional
(\ref{31}) is a function of values of $Y$ at finite number of time
moments, thus the mean value in (\ref{32}) is well defined  as the
integral over the family $\{P_{st_1\dots t_k}(x,\cdot),
t_1,\dots,t_k>s, k\in \NN\}$ of conditional finite-dimensional
distributions  of the process $Y$.

The main result of this chapter is given in the following
theorem.

\begin{thm}\label{t31}
 Assume that there exist the sequence $X_n$ that provides Markov approximation
 for the  homogeneous Markov process $X$ and the sequence
  $\{\phi_n\equiv \phi_n(X_n)\}$ of the functionals of the type (\ref{31}).
 Let the following conditions hold true:
\begin{enumerate}
\item The functions $F_{n}(\cdot)$ are bounded  and uniformly tend
to zero:
$$\delta(F_{n})\eqdef \sup\{F_{n}(x_{1},\dots,x_{L})|{x_{1},\dots,x_{L}\in \Xf}\}
\rightarrow0,\quad n\rightarrow\infty.$$ \item There exists a
function $f$, that appears to be a characteristics (in a sense of
Chapter 6 \cite{dynkin}) of some  $W$-functional $\phi=\phi(X)$ of
the limiting  Markov process  $X$, such that, for each  $T$,
$$\underset{s={i\over n}, t\in (s,T)}{\sup}\norm{f^{s,t}_n(\cdot)-f^{s,t}(\cdot)}\rightarrow0,
\quad n\rightarrow\infty,$$ where $\norm{g(\cdot)}\equiv
\underset{x\in \Xf}\sup|g(x)|$.

\item  The limiting function  $f$ is uniformly continuous with
respect to variable $x$, that is, for arbitrary  $T$
$$\underset{0\leq s\leq t<T}{\sup}\abs{f^{s,t}(x')-f^{s,t}(x'')}
\rightarrow0,\quad \abs{x'-x''}\rightarrow0.$$
\end{enumerate}
Then
$$\psi_n(X_n)\Rightarrow\phi(X)\equiv \{\phi^{s,t}(X), (s,t)\in \TT\},$$
where  $\psi_n$ are the random broken lines   corresponding to the
functionals $\phi_n$ and convergence  is understood in a sense of
$C(\TT,\Re^+)$.
\end{thm}

\begin{rem} Conditions 1,2 are analogous to those of the Dynkin's theorem:
condition 2 is exactly the condition for the characteristics to
converge, condition 1 corresponds to the assumption that the
prelimit functionals are $W$-functionals. In the present
situation, of course, we can not say that $\phi_n$ are
$W$-functionals, particulary, $\phi_n$ are not continuous with
respect to temporary variable. Condition 1 means exactly that the
values of jumps are negligible while $n\to \infty$. Condition 3,
though not very restrictive, is specific, and is caused by
necessity to consider functionals, set over different processes.
\end{rem}
\begin{rem}\label{r32} If $X_n\Rightarrow X$ in  $\CC$ or in
$\DD$ (this condition is not provided by the conditions of the
Theorem,  see Remark \ref{r13}), then, as one can easily see from
the proof,  $(X_n,\psi_n(X_n))\Rightarrow (X,\phi(X))$ in
$\CC\times C(\TT,\Re^+)$ or in $\DD\times C(\TT,\Re^+)$,
respectively.
 \end{rem}

Note that the result of the theorem also holds for the Markov
process $X$ that is not homogeneous w.r.t. time variable; the
claim for the limit Markov process to be homogeneous is imposed in
order to shorten the notation only. This remark concerns also the
most of the results stated below.

{\it Proof of the theorem.} The general scheme of the proof is
close to the  one, proposed in \cite{andr_kulik} in order to prove
the analogue of the Dynkin's theorem for the family of functionals
of a single Markov process, for which the properties of
additivity, continuity and homogeneity may fail, but the
violations become negligible while $n\to \infty$.

First let us show that the  finite-dimensional distributions of
$\phi_n$ converge to the corresponding distributions of $\phi$.
Let the constants $\gamma,T$ be fixed and
$\widehat{X}_n,\widehat{X}^n$ be  processes satisfying conditions
(i)-(iii) of Definition \ref{d11} with these constants. For these
processes, one can consider the functionals $\phi_n(\hat
X_n),\phi(\hat X^n)$; obviously, their distributions and
characteristics coincide with those for $\phi_n(X_n),\phi(X)$. In
order to shorten notation, we denote further $\phi_n=\phi_n(\hat
X_n),\phi=\phi(\hat
X^n),K=K(\gamma,T),\Fc_t=\hat\Fc_t^n\equiv\sigma(\hat X_n(s),\hat
X^n(s),s\leq t).$

It follows from the condition (iii) and the definition of
characteristics that, for arbitrary $t\in \left({iK\over
n},T\right]$,
\be\label{33}E\left[\phi^{\frac{Ki}{n},t}|\mathcal{F}_{\frac{Ki}{n}}
\right]=f^{\frac{Ki}{n},t}\brakes{\widehat{X}^n\brakes{\frac{K
i}{n}}},\quad
E\left[\phi^{\frac{Ki}{n},t}_n|\mathcal{F}_{\frac{Ki}{n}}\right]=
f_n^{\frac{Ki}{n},t}\brakes{\widehat{X}_n\brakes{\frac{K
i}{n}}}\ee almost surely.

\begin{lem}\label{l31} For $0\leq s\leq t\leq T$, the following estimate holds:
$$\lims E\brakes{\phi^{s,t}_{n}(\widehat{X}_n)-\phi^{s,t}(\widehat{X})}^2\leq
4\norm{f^{0,T}} G(f,\gamma,T)+4\sqrt{2\gamma}\norm{f^{0,T}}^2,$$
where $G(f,\gamma,T)=\underset{0\leq s\leq t\leq T,
\abs{x'-x''}<\gamma}{\sup}\abs{f^{s,t}(x')-f^{s,t}({x''})}$.
\end{lem}

{\it Proof.} We will prove the statement of lemma for $s=0, t=T$;
in general case the proof is exactly the same.
  Consider the partition of the axis $\ax$ by points of the type
$\frac{Ki}{n}, i\in \NN$.  Denote $M_n=[\frac{nT}{K}]+1$,
$$\Delta_i^n\eqdef \phi^{\frac{(i-1)K}{n},(\frac{iK}{n})\bigwedge
T}_n,\quad\widetilde{\Delta}_i^n\eqdef
\phi^{\frac{(i-1)K}{n},\frac{iK}{n}\bigwedge T},\quad
i=\overline{1,M_n}.$$

We have that
$$\brakes{\phi^{0,T}_n-\phi^{0,T}}^2=\brakes{\sum_{i=1}^{M_n}\Delta_i^n-\widetilde{\Delta}_i^n}^2=
\brakes{\sum_{i=1}^{M_n}\Delta_i^n}^2+\brakes{\sum_{i=1}^{M_n}\widetilde{\Delta}_i^n}^2-
2\sum_{i=1}^{M_n}\sum_{j=1}^{M_n}\Delta_i^n\widetilde{\Delta}_j^n=\Sigma_1^n+2\Sigma_2^n,$$
where
$$\Sigma_1^n\eqdef \sum_{i=1}^{M_n}(\Delta_i^n)^2+\sum_{i=1}^{M}(\widetilde{\Delta}_i^n)^2-2\sum_{i=1}^{M}\Delta_i^n
\widetilde{\Delta}_i^n,$$
$$\Sigma_2^n\eqdef\sqrbrakes{\sum_{1\leq i<l\leq M_n}\Delta_i^n\Delta_l^n-\sum_{1\leq i<j\leq M_n}\Delta_i^n
\widetilde{\Delta}_j^n}+\sqrbrakes{\sum_{1\leq j<k\leq
M_n}\widetilde{\Delta}_j^n\widetilde{\Delta}_k^n- \sum_{1\leq
j<i\leq M_n}\Delta_i^n\widetilde{\Delta}_j^n}.$$
Let us estimate the expectations $\Sigma_1^n,\Sigma_2^n$
separately. Since the increments
$\Delta_i^n,\widetilde{\Delta}_i^n$ are non-negative, the first
sum can be estimated  by the sum of the  first two terms:
\be\label{34}\Sigma_1^n\leq\sum_{i=1}^{M_n}(\Delta_i^n)^2+\sum_{i=1}^{M_n}(\widetilde{\Delta}_i^n)^2.\ee
The expectation of the first term in (\ref{34}) can be estimated
via the definition of  $\phi_n$:
$$E\sum_{i=1}^{M_n}(\Delta_i^n)^2\leq
E\brakes{\sup\limits_{i=1,M_n}\Delta_i}\sum_{i=1}^{M_n}\Delta_i^n\leq
K\delta_n
 f^{0,T}_n\brakes{\hat X_n(0)}\leq K\delta_n\norm{f^{0,T}_n}\to 0, \quad n\to +\infty,
$$
 where $\delta_n\equiv\delta(F_{n})$.
Convergence to zero of the expectation of the second term in
(\ref{34}) is provided by the arguments, analogous to those used
in \cite{dynkin} Chapter 6:  on the one hand, by the continuity of
functional $\phi$,
$\sum_{i=1}^{M_n}(\widetilde{\Delta}_i^n)^2\to0$ by probability;
on the  other hand, $\sum_{i=1}^{M_n}(\widetilde{\Delta}_i^n)^2$
is dominated by the variable $(\phi^{0,T})^2$; the expectation of
this variable, due to Lemma 6.4 \cite{dynkin}, does not exceed
$2\norm{f^{0,T}}^2<\infty$. Therefore,
$E\sum_{i=1}^{M_n}(\widetilde{\Delta}_i^n)^2\to0$ due to the
Lebesgue theorem on dominated convergence. Hence, $\lims
E\Sigma^n_1\leq 0$.


   The expectation of $\Sigma_2^n$ is equal
$$E\Sigma_2^n=E\sqrbrakes{\sum_{1\leq i<l\leq M_n}\Delta_i^n\Delta_l^n-\sum_{1\leq i<j\leq M_n}\Delta_i^n\widetilde{\Delta}_j^n}+
E\sqrbrakes{\sum_{1\leq j<k\leq
M_n}\widetilde{\Delta}_j^n\widetilde{\Delta}_k^n-\sum_{1\leq
j<i\leq M_n}\Delta_i^n \widetilde{\Delta}_j^n}=$$
\be\label{35}=E\sum_{i=1}^{M_n-1}\Delta_i^n\sqrbrakes{\phi^{\frac{Ki}{n},T}_n-\phi^{\frac{Ki}{n},T}}-
E\sum_{i=1}^{M_n-1}\widetilde{\Delta}_i^n\sqrbrakes{\phi^{\frac{Ki}{n},T}_n
-\phi^{\frac{Ki}{n},T}}.\ee We estimate the second term in
(\ref{35}), using property (\ref{33}). Since  $\tilde \Delta_i^n$
is measurable w.r.t. $\mathcal{F}_{\frac{Ki}{n}}$, the following
estimate holds:
$$-E\sum_{i=1}^{M_n-1}\widetilde{\Delta}_i^n\left[\phi^{\frac{Ki}{n},T}_n
-\phi^{\frac{Ki}{n},T}\right]=-E\sum_{i=1}^{M_n-1}\widetilde{\Delta}_i^nE\left[\brakes{\phi^{\frac{Ki}{n},T}_n-\phi^{\frac{Ki}{n},T}}|\mathcal{F}_{\frac{Ki}{n}}\right]\leq$$
$$\leq
E\sum_{i=1}^{M_n-1}\widetilde{\Delta}_i^n\brakes{f^{\frac{Ki}{n},T}_n\brakes{\hat
X_n\brakes{{Ki\over n}}}-f^{\frac{Ki}{n},T}\brakes{\hat
X^n\brakes{{Ki\over n}}}}\leq
$$
$$\leq
E\sum_{i=1}^{M_n-1}\widetilde{\Delta}_i^n\abs{f_n^{\frac{Ki}{n},T}-f^{\frac{Ki}{n},T}}+E
\sum_{i=1}^{M_n-1}\widetilde{\Delta}_i^n\abs{f^{\frac{Ki}{n},T}\brakes{\hat
X_n\brakes{{Ki\over n}}}-f^{\frac{Ki}{n},T}\brakes{\hat
X^n\brakes{{Ki\over n}}}}\leq
$$
\be\label{36}\leq \|f^{0,T}\| \underset{s={i\over n}, t\in
(s,T)}{\sup}\norm{f^{s,t}_n(\cdot)-f^{s,t}(\cdot)}+E
\sum_{i=1}^{M_n-1}\widetilde{\Delta}_i^n\abs{f^{\frac{Ki}{n},T}\brakes{\hat
X_n\brakes{{Ki\over n}}}-f^{\frac{Ki}{n},T}\brakes{\hat
X^n\brakes{{Ki\over n}}}}
 \ee
(in the last inequality, we have used that
$\sum_{i=1}^{M_n-1}\widetilde{\Delta}_i^n\leq \phi^{0,T}$ and
$E\phi^{0,T}\leq \|f^{0,T}\|$). The first term in (\ref{36}) tends
to zero. In order to estimate the second term, we put
$\Omega_{\gamma,T}=\left\{\sup\limits_{i\leq {Tn\over K}}
\rho\left(\hat X_n\left({i K\over n}\right),\hat X^n\left({i
K\over n}\right)\right)>\gamma\right\}$ (recall that
$P(\Omega_{\gamma,T})<\gamma$ due to the claim (iii) of Definition
\ref{d11}). We have
$$E
\sum_{i=1}^{M_n-1}\widetilde{\Delta}_i^n\abs{f^{\frac{Ki}{n},T}\brakes{\hat
X_n\brakes{{Ki\over n}}}-f^{\frac{Ki}{n},T}\brakes{\hat
X^n\brakes{{Ki\over n}}}}\leq
E\phi^{0,T}G(f,\gamma,T)\1_{\Omega\backslash\Omega_{\gamma,T}}+
$$
\be\label{37}
+E\sum_{i=1}^{M_n-1}\widetilde{\Delta}_i^n\abs{f^{\frac{Ki}{n},T}\brakes{\hat
X_n\brakes{{Ki\over n}}}-f^{\frac{Ki}{n},T}\brakes{\hat
X^n\brakes{{Ki\over n}}}}\1_{\Omega_{\gamma,T}}.\ee The first term
in (\ref{37}) can be estimated by $\|f^{0,T}\| G(f,\gamma,T)$. The
second term is estimated by Cauchy inequality:
$$
E\sum_{i=1}^{M_n-1}\widetilde{\Delta}_i^n\abs{f^{\frac{Ki}{n},T}\brakes{\hat
X_n\brakes{{Ki\over n}}}-f^{\frac{Ki}{n},T}\brakes{\hat
X^n\brakes{{Ki\over n}}}}\1_{\Omega_{\gamma,T}}\leq
$$
$$
\leq\norm{f^{0,T}}E\phi^{0,T}\1_{\Omega_{\gamma,T}}\leq
\norm{f^{0,T}}\sqrbrakes{E(\phi^{0,T})^2}^{1\over
2}\sqrbrakes{P(\Omega_{\gamma,T})}^{1\over
2}\leq\norm{f^{0,T}}^2\sqrt{2\gamma}
$$
(here, the Lemma 6.4 \cite{dynkin} was applied). Summing up the
above relations, we deduce that \be\label{38}
\lims\left\{-E\sum_{i=1}^{M_n-1}\widetilde{\Delta}_i^n\sqrbrakes{\phi^{\frac{Ki}{n},T}_n
-\phi^{\frac{Ki}{n},T}}\right\}\leq \norm{f^{0,T}}
G(f,\gamma,T)+\norm{f^{0,T}}^2\sqrt{2\gamma}.\ee

Now, let us proceed with the estimation of the first item in
(\ref{35}). Straightforward use of the property (\ref{33}) is
impossible here, since the variable $\Delta_i^n$ is a functional
of values of the process $\hat X_n$ at the points ${Ki\over
n},{Ki+1\over n},\dots{Ki+L\over n}$, that is, it is not
measurable with respect to $\mathcal{F}_{\frac{Ki}{n}}$. Without
loss of generality, one can  assume that $K\geq L$ (otherwise one
can make the same procedure with the constant $K$ replaced by
 $K\cdot L$). Then the variable $\Delta_i^n$ is measurable with respect
 to $\mathcal{F}_{\frac{K(i+1)}{n}}$. The functionals
$\phi_n,\phi$ are additive at points of the type ${j\over n}$.
Applying (\ref{33}) and condition 1 of the Theorem, we obtain the
following relation
$$
E\sum_{i=1}^{M_n-1}\Delta_i^n\sqrbrakes{\phi^{\frac{Ki}{n},T}_n-\phi^{\frac{Ki}{n},T}}=
E\sum_{i=1}^{M_n-1}\Delta_i^n\sqrbrakes{\phi^{\frac{Ki}{n},\frac{K(i+1)}{n}}_n-\phi^{\frac{Ki}{n},\frac{K(i+1)}{n}}}+
$$
$$
+E\sum_{i=1}^{M_n-1}{\Delta}_i^n\sqrbrakes{f^{\frac{K(i+1)}{n},T}\brakes{\hat
X_n\brakes{{K(i+1)\over n}}}-f^{\frac{K(i+1)}{n},T}\brakes{\hat
X^n\brakes{{K(i+1)\over n}}}}\leq
$$
 \be\label{39}\leq K\delta_n\abs{f^{0,T}_n}+E\sum_{i=1}^{M_n-1}{\Delta}_i^n\sqrbrakes{f^{\frac{K(i+1)}{n},T}\brakes{\hat
X_n\brakes{{K(i+1)\over n}}}-f^{\frac{K(i+1)}{n},T}\brakes{\hat
X^n\brakes{{K(i+1)\over n}}}}. \ee The first term in (\ref{39})
tends to zero. The second term in (\ref{39}) is estimated in  the
same way with the second term in (\ref{35}), with one necessary
change. We cannot apply Lemma 6.4 \cite{dynkin} in order to
estimate the second moment $\phi_n^{0,T}$, therefore this estimate
must be obtained separately. This can be done in a following way:
$$E(\phi_n^{0,T})^2=E\sum_{i=1}^{M_n}(\Delta_i^n)^2+2E\sum_{1\leq
i< j\leq M_n}^{M_n}\Delta_i^n\Delta_j^n=
E\sum_{i=1}^{M_n}(\Delta_i^n)^2+2E\sum_{1\leq i\leq
M_n}^{M_n}\Delta_i^n\phi_n^{\frac{iK}{n},T}=$$
$$=E\sum_{i=1}^{M_n}(\Delta_i^n)^2+2E\sum_{1\leq i\leq M_n}^{M_n}\Delta_i^n[\phi_n^{iK/n,(i+1)K/n}+\phi_n^{(i+1)K/n,T}]
\leq$$
$$=E\sum_{i=1}^{M_n}(\Delta_i^n)^2+2K\delta_nE\sum_{1\leq i\leq M_n}^{M_n}\Delta_i^n+
2E\sum_{1\leq i\leq M_n}^{M_n}\Delta_i^n\norm{f_n^{0,T}}\leq$$
\be\label{391}\leq\left\{(2K+1)\delta_n+2\|f_n^{0,T}\|\right\}E\phi_n^{0,T}\leq
(2K+1)\delta_n\norm{f_n^{0,T}}+2\norm{f_n^{0,T}}^2,\ee
 all transitions here are analogous to those  given above,
 and thus are not discussed in details.
 Repeating literally the estimates for the second term in (\ref{35}),
 we obtain the estimate \be\label{310}
\lims
E\sum_{i=1}^{M_n-1}\Delta_i^n\sqrbrakes{\phi^{\frac{Ki}{n},T}_n-\phi^{\frac{Ki}{n},T}}\leq
\norm{f^{0,T}} G(f,\gamma,T)+\norm{f^{0,T}}^2\sqrt{2\gamma}. \ee
It follows from (\ref{38}),(\ref{310})  that
$\lims\sqrbrakes{2\Sigma_2^n}\leq 4\norm{f^{0,T}}
G(f,\gamma,T)+4\sqrt{2\gamma}\norm{f^{0,T}}^2$. This, combined
with the estimate $\lims\sqrbrakes{\Sigma_1^n}\leq 0$ proved
before, provides the needed statement. The lemma is proved.

Now, we can complete the proof of the convergence of
finite-dimensional distributions of $\phi_n$  to those of $\phi$.
In order to shorten notation we consider the one-dimensional
distributions only; in general case considerations are completely
the same.

 Take arbitrary $s,t, s<t$. In order to prove weak convergence
 $\phi_n^{s,t}(X_n)$ to $\phi^{s,t}(X)$, it is sufficient to show
 that, for arbitrary bounded Lipschitz function $g$,
\be\label{311} \lims
\abs{Eg(\phi_n^{s,t}(X_n))-Eg(\phi^{s,t}(X))}=0. \ee Let $g$ be
fixed, consider a pair of processes $\hat X_n,\hat X^n$,
corresponding (in a sence of Definition \ref{d11}) to $T=t$ and
given positive $\gamma$. By construction, $\phi_n^{s,t}(X_n)\eqd
\phi_n^{s,t}(\hat X_n),\phi^{s,t}(X)\eqd \phi^{s,t}(\hat X^n)$.
Applying Lemma \ref{l31}, we obtain that
$$ \lims \abs{Eg(\phi_n^{s,t}(X_n))-Eg(\phi^{s,t}(X))}\leq \lims
E\abs{g(\phi_n^{s,t}(\hat X_n))-\phi^{s,t}(\hat X^n)}\leq
$$
$$
\leq \mathrm{Lip}(g)\,\lims E\abs{\phi_n^{s,t}\brakes{\hat
X_n}-\phi^{s,t}\brakes{\hat X^n}}\leq
2\mathrm{Lip}(g)\sqrt{\norm{f^{0,t}}
G(f,\gamma,t)+\sqrt{2\gamma}\norm{f^{0,t}}^2},
$$
here $\mathrm{Lip}(g)$ denotes the Lipshits constant for $g$.
Condition 3 of the Theorem provides  that $G(f,\gamma,t)\to
0,\gamma\to 0+$. Therefore, since  $\gamma>0$ is arbitrary,
(\ref{311}) follows from the estimate given above.

Since $\sup_{s,t}|\psi_n^{s,t}-\phi_n^{s,t}|\leq \delta_n\to 0$,
the finite-dimensional distributions of $\phi_n$ converge to
corresponding distributions of $\phi$. Thus, the only thing left
to show in order to prove the Theorem, is that the family of
distributions of $\psi_n$
 is dense in $C(\TT,\Re^+)$. The values of the functions $\psi_n$ at the point $s,t$
differ from the values at the closest knots of partition
$s_*,t_*\in {1\over n}\ZZ_+$ at most on $\delta_n$, and $\psi_n$
are monotone as the functions of the time variables. Hence, in
order to prove the required statement, it is sufficient to show
that, for arbitrary sequence of partitions
$\left\{S_n=\{s_0^n=0<s_1^n<\dots<s_k^n<\dots\}\subset {1\over
n}\ZZ_+,n\in\NN\right\}$ with  $\sigma_n\equiv
\max_k(s_k^n-s_{k-1}^n)\to 0, n\to +\infty$ and arbitrary $T\in
\ax$,
$$
E\sum_{k:s_k\leq T}\left[\psi_n^{s_{k-1}^n,s_k^n}\right]^2\to
0,\quad n\to +\infty.
$$
Set $\gamma_{n,T}=\sup_{0<t-s<\sigma_n, t<T}\|f_n^{s,t}\|$, note,
that $\gamma_{n,T}\to 0,n\to +\infty$ due to continuity of the
limit characteristics $f$ and uniform convergence of
$f_n\rightrightarrows f$. In the same way with (\ref{391}) we
obtain the estimate \be\label{312}
E\left[\phi_n^{s_{k-1}^n,s_k^n}\right]^2\leq
\left\{(2K+1)\delta_n+2\gamma_{n,T}\right\}E\phi_n^{s_{k-1}^n,s_k^n}.
\ee  Summing up  the estimates (\ref{312}) w.r.t. $k$ (recall that
$\phi_n^{s,t}=\psi_n^{s,t}$ when $s,t\in {1\over n}\ZZ_+$), we
obtain
$$
E\sum_{k:s_k\leq T}\left[\psi_n^{s_{k-1}^n,s_k^n}\right]^2\leq
\left\{(2K+1)\delta_n+2\gamma_{n,T}\right\}\|f_n^{0,T}\|\to
0,\quad n\to +\infty,$$ what was to be proved. The theorem is
proved.

Let us make one remark. For the random walks, the Skorokhod's
method is well known, allowing one to reduce the investigation of
the sums of the type (\ref{31}) to the case $L=1$. This method can
be applied in the context of current paper, also. Namely, the
reasoning, similar to the one used in the proof of Theorem 1,
Chapter 5.3 \cite{skor_slob}, provides the following result (the
proof is omitted).

\begin{prop}\label{c31} Let the sequence of functionals $\{\phi_n=\phi_n(X_n)\}$ of the type
(\ref{31}) be given, and, for every $n$, the process $X_n$
possesses the Markov property at the time moments ${i\over n},
i\in \ZZ_+$. Consider the functionals
$$\chi_{n}^{s,t}(X_n)\eqdef \sum_{k:s\leq
k/n<t}\Psi_{n,k}\brakes{X_n\brakes{\frac{k}{n}}},\quad 0\leq s<t,
$$
where
$$
\Psi_{n,k}(x)\equiv E\left[F_{n}
\brakes{x,X_n\brakes{\frac{k+1}{n}},\dots,X_n\brakes{\frac{k+L-1}{n}}}\Bigl|
X_n\brakes{\frac{k}{n}}=x\right],\quad x\in \Xf.
$$
Let functions $F_{n}(\cdot)$ be non-negative and satisfy condition
1 of Theorem \ref{t31}, then the functionals $\phi_n^{s,t}$ have a
limit distribution if and only if the functionals $\chi_n^{s,t}$
have a limit distribution, and the limit  distributions of the
functionals $\phi_n^{s,t},\chi_n^{s,t}$ are equal as soon as they
exist.
\end{prop}

It is worth to note that the Proposition \ref{c31} does not lead
to simplification of the initial problem in the context of current
paper. The number of values of process $X_n$, contained in a one
summand for the functional $\phi_n$ (that is, number $L$), is not
involved significantly into the proof of the main theorem. We will
see later that the  main problem in the application of the Theorem
consists in verification of the condition 2 of uniform convergence
of characteristics; the characteristics of the functionals
$\phi_n$ and $\chi_n$, obviously, coincide.

In the following two chapters, the examples of application of
Theorem \ref{t31} are given.

\section{The local time of a random walk at a point.}

Let the processes $X_n$ be constructed w.r.t. one-dimensional
random walk that belongs to the normal domain of attraction of an
$\alpha$-stable law, $\alpha\in(1,2]$ (see Examples 1,2). We
assume the centering sequence $a_n$ to be equal to zero, and set
the random broken lines $X_n$ by equality (\ref{22}).

Consider, for arbitrary $z_*\in\Re$, the  functionals
$\phi_n=\phi_n(X_n)$ of the type (\ref{31}) with
$L=2,F_{n}(x,y)={1\over n}{1\over
|y-x|}\left[\1_{(x-z_*)(y-z_*)<0}+{1\over
2}(\1_{x\not=z_*,y=z_*}+\1_{x=z_*,y\not=z_*})\right]$. For every
$s<t, s,t\in \{{j\over n},j\in \ZZ_+\}$, with probability 1 the
following equality takes place \be\label{400}
\phi_n^{s,t}(X_n)=\lim_{\eps\to 0+}{1\over
2\eps}\int_s^t\1_{X_n(r)\in(z_*-\eps,z_*+\eps)\backslash\{z_*\}}\,
dr,\quad 0\leq s<t. \ee Therefore the functionals $\phi_n$ can be
naturally interpreted as the \emph{censored} local times for the
broken lines $X_n$ at the point $z_*$ (the censoring operation
consists in removing horizontal parts of the broken lines).
Theorem 3.1 allows one to obtain the following limit result.

\begin{prop}\label{c41} Let the
distribution of the jump $\xi_1$ of the random walk be
concentrated on ${\ZZ}$ and aperiodic. Then the conditions of
Theorem \ref{t31} hold true and $\phi_n^{s,t}(X_n)$ converge by
distribution to $\phi^{s,t}(X)=P(\xi_1\not=0)\cdot
L^{s,t}(X,z_*)$, where $L(X,z_*)$ is the local time of the limit
$\alpha$-stable process $X$ at the point $z_*$.
\end{prop}

{\it Proof.} The condition for $X_n$ to provide Markov
approximation  for $X$ holds true (see Example 2). Condition 1 of
the  Theorem holds with $\delta_n=2n^{{1\over \alpha}-1}$ since
either the increment of the process $X_n$ in the neighboring knots
is equal to zero or the absolute value of this increment  is not
less then $n^{-{1\over \alpha}}$. Let us show that the
characteristics of functionals $\phi_n$ converge uniformly to the
function \be\label{40}
f^{s,t}(x)=P(\xi_1\not=0)\int_0^{t-s}p_r(z_*-x)\,dr, \ee where
$p_r(\cdot)$ is the density of distribution $X(r)$ under condition
$X(0)=0$; this provides conditions 2,3 of the Theorem.

In order to shorten notation we take $z_*=0$. Denote
$P^k_i=P(S_k=i), P_j=P^1_j=P(\xi_1=j), i,j\in \ZZ$. We have that
$$
f_n^{s,t}(x)=n^{{1\over \alpha}-1} \sum_{s\leq {k\over
n}<t}\left[\sum_{j\not=0}{P_j\over |j|}\left(\sum_{i\in
(xn^{1\over \alpha}-j,xn^{1\over \alpha})}P^k_i+{1\over
2}\1_{xn^{1\over \alpha}\in \ZZ}\brakes{P_{xn^{1\over
\alpha}}^k+P_{xn^{1\over \alpha}-j}^k}\right)\right],
$$
notation $i\in(a,b)$ in the case $a>b$  means that $b<i<a$. Using
the appropriate version of the Gnedenko's local limit theorem (see
\cite{ibr_linnik}, Theorem 4.2.1), one can write \be\label{41}
\eps_k\equiv \sup_{i\in \ZZ}\left|k^{1\over \alpha}P_i^k -
p_1\brakes{{i\over k^{1\over \alpha}}}\right|\to 0,\quad k\to
+\infty. \ee Hence
$$
f_n^{s,t}(x)={1\over n}\sum_{s\leq {k\over
n}<t}\left[\sum_{j\not=0}{P_j\over |j|}\left(\sum_{i\in
(xn^{1\over \alpha}-j,xn^{1\over \alpha})}\brakes{{n\over
k}}^{1\over \alpha}p_1\brakes{i\over k^{1\over \alpha}}+
\right.\right.
$$
\be\label{42} +\left.\left. {n^{1\over \alpha}\over 2k^{1\over
\alpha}}\1_{xn^{1\over \alpha}\in \ZZ}\left\{p_1\brakes{xn^{1\over
\alpha}\over k^{1\over \alpha}}+p_1\brakes{xn^{1\over
\alpha}-j\over k^{1\over \alpha}}\right\}\right)\right]+\Xi_n(x),
\ee where \be\label{4_estimate} |\Xi_n(x)|\leq {1\over
n}\sum_{k=1}^{[nt]}\brakes{n\over k}^{1\over \alpha}\eps_k, \ee
and  $\Xi_n\rightrightarrows 0,n\to +\infty$ via the Toeplitz's
theorem.

The density $p_1$ is uniformly continuous over $\Re$, hence, using
the same arguments, one can show that, up to a summand that
uniformly converges to zero, the value of $f^{s,t}_n(x)$ equals
$$
{1\over n}\sum_{s\leq {k\over n}<t}\left[\sum_{j\not=0}{P_j\over
|j|}\left(\sum_{i\in (xn^{1\over \alpha}-j,xn^{1\over
\alpha})}\brakes{{n\over k}}^{1\over \alpha}p_1\brakes{xn^{1\over
\alpha}\over k^{1\over \alpha}}+  \brakes{{n\over k}}^{1\over
\alpha}\1_{xn^{1\over \alpha}\in \ZZ}p_1\brakes{xn^{1\over
\alpha}\over k^{1\over \alpha}}\right)\right]=
$$
\be\label{43} ={P(\xi_1\not=0)\over n}\sum_{s\leq {k\over n}<t}
\brakes{{n\over k}}^{1\over \alpha}p_1\brakes{xn^{1\over
\alpha}\over k^{1\over \alpha}}={P(\xi_1\not=0)\over n}\sum_{s\leq
{k\over n}<t}p_{k\over n}(x); \ee  in the latter equality, we have
used  that the process $X$ is self-similar, that is,
$p_r(x)=r^{-{1\over \alpha}}p_1(r^{-{1\over \alpha}}x), r>0$. The
sum in the right hand part of (\ref{43}) is exactly the integral
sum for the integral in the right hand part of (\ref{40}), the
functions $\{p_r(\cdot), r\geq r_0\}$ are uniformly continuous for
arbitrary $r_0>0$ and $\sup_xp_r(x)\leq Cr^{-{1\over \alpha}}$.
This immediately provides the required uniform convergence of
$f_n$ to $f$. The proposition is proved.

The similar result can be proved for $\xi_k$ with  non-lattice
distribution, for which there exists a bounded  distribution
density of $S_{n_0}$ for some $n_0$ (the  proof is omitted).

The result of Proposition  \ref{c41} and its analog for
non-lattice random walks is not essentially new;  one can obtain
it applying either Proposition \ref{c31} and the technique,
exposed in \S\S III.2, III.3 \cite{bor_ibrag}, or the reasonings,
similar to those used in the proof of Theorem 3
\cite{portenko_68}. Our reason to give this example consists, on
the one hand, in describing the way of application of Theorem
\ref{t31}  in a simple situation where an appropriate local limit
theorem is available, and on the other hand, in emphasizing the
following interesting fact, that is not reflected in a literature
available for us. For a "good"\phantom{} random walks (lattice or
essentially non-lattice), their local times at the point, defined
by the natural equality (\ref{400}), converge by distribution
exactly to the local time of the limit process at the same point,
as soon as the broken lines corresponding to the random walk does
not contain horizontal sections.

\section{Difference approximations of diffusion processes.}
Consider the sequence $\{Z_n\}$ of difference approximations of
diffusion process $Z$ (see Example \ref{e3}, equalities
(\ref{23}),(\ref{24})). The sequence $\{Z_n\}$ provides Markov
approximation for $Z$, that allows one to apply Theorem \ref{t31}
while considering the question on the  limit behavior of the
functionals of type (\ref{31}) for $\{Z_n\}$.

One of possible way to proceed here is to apply the estimates
based on an appropriate local limit theorem, like it was made in
the previous chapter. In order to make this paper reasonably
short, we do not give the detailed  exposition of this subject
here (see the separate paper \cite{kulik_prep}). In this chapter,
we give a simple corollary of Theorem \ref{t31}, that provides
invariance principle for certain "canonic"\phantom{} additive
functionals, that are related to the Doob's decomposition of
$|Z_n(\cdot)|$.

Let us consider the objects introduced in Example \ref{e3} with
$m=d=1$ and $a,b,\{\xi_n\}$ satisfying conditions introduced
there. Put \be\label{431} \phi_{n}^{s,t}(Z_n)\equiv\!
\sum_{k\in(sn,tn]}\!|Z_n\brakes{{k\over
n}}|\cdot\left[2\1_{Z_n\brakes{{k-1\over n}}Z_n\brakes{{k\over
n}}< 0}+\1_{Z_n\brakes{{k-1\over n}}=0}\right], \ee $\psi_n$ are
corresponding broken lines.

\begin{prop}\label{c42} The processes $\psi_n$ converge by distribution in $C(\TT,\Re)$
to the local time
$$
\phi^{s,t}\equiv \lim_{\eps\to 0+}{1\over
2\eps}\int_s^t\1_{|Z(r)|<\eps}b^2(Z(r))\,dr
$$
 of the diffusion process $Z$ at the point $0$.
\end{prop}

{\it Proof.} Since the diffusion coefficient is non-degenerate,
$Z$ possesses continuous transition density $p_t(x,y)$  and the
standard estimate $\sup_xp_t(x,y)\leq {C(y)\over \sqrt t}$ holds
true. This implies existence of the local time of $Z$ at the point
$0$. This local time is a $W$-functional with the characteristics
$ f^{0,t}(x)=b^2(0)\int_0^tp_{s}(x,0)\,ds, $ that is, condition 3
of Theorem \ref{t31} holds. Straightforward calculations prove the
equality \be\label{44}
|Z_n(t)|-|Z_n(s)|=\phi_{n}^{0,t}(Z_n)+\sum_{k=ns}^{[nt]-1}\left[a\brakes{Z_n\brakes{{k\over
n}}}{1\over n}+b\brakes{Z_n\brakes{{k\over n}}}\Delta
X_n\brakes{{k\over n}}\right]\sign \brakes{Z_n\brakes{{k\over
n}}}, \ee where $s\in{1\over n}\ZZ_+,\sign(0)=0$. This provides
that
$$
f_n^{s,t}(x)=E\left[|Z_n(t)||Z_n(s)=x\right]-|x|- {1\over
n}E\left[\sum_{k=0}^{[nt]-1}a\brakes{Z_n\brakes{{k\over n}}}\sign
(Z_n\brakes{{k\over n}}\Bigl|Z_n(s)=x\right].
$$
Processes $Z_n$ converge weakly to $Z$, function $a(x)\sign(x)$
has unique jump at point $x=0$ and  $P(Z(r)=0)=0$ for every $r>0$.
Hence the standard reasonings provide that (we omit the details)
\be\label{45} f_n^{s,t}(x)\mathop{\rightrightarrows}\limits_x
E\left[|Z(t)||Z(s)=x\right]-|x| -
E\left[\int_s^ta(Z_r)\sign(Z_r)\,dr\Bigl|Z(s)=x\right]. \ee This
proves condition 2 of Theorem \ref{t31}, since the right hand side
of (\ref{45}) is exactly the characteristics of the local time
$\phi$ due to Ito-Tanaka formula.

In order to provide  condition 1, let us, for a while, suppose
additionally that the coefficients $a,b$ are bounded. We apply the
standard "cutting" procedure: on each step of approximation,
together with the process $Z_n$, we consider the process $\tilde
Z_n$, constructed by the same scheme from a sequence of i.i.d.r.v.
$\{\tilde\xi_n\}$, satisfying conditions $\|\tilde \xi_n\|\leq
n^{{1\over 2+{\delta\over 2}}}$ and $\xi_n=\tilde \xi_n$ for
$\|\xi_n\|\leq n^{{1\over 2+{\delta\over 2}}}$. For such $\tilde
Z_n$, condition 1 of theorem holds with
$$\delta(F_{n})\leq n^{-1}\max_x|a(x)|+n^{{1\over 2+{\delta\over
2}}-{1\over 2}}\max_x|b(x)|,$$ and the other conditions of theorem
for $\tilde Z_n$ remain to hold true. This proves the statement of
Proposition \ref{c42} for $\{\tilde Z_n\}$. On the other hand, for
arbitrary $T\in\ax$
$$
P\brakes{Z_n|_{[0,T]}\not=\tilde
Z_n|_{[0,T]}}=O\brakes{n^{1-{2+\delta\over 2+{\delta\over
2}}}}=o(1), \quad n\to +\infty,
$$
and therefore the statement of Proposition \ref{c42} holds true
for $\{Z_n\}$. At last, the additional assumption that  the
coefficients $a,b$ are bounded, can be removed via a standard
localization procedure. The proposition is proved.

\begin{rem} Let $a=0,b=1,P(\xi_k=\pm 1)={1\over 2}$ (that is, $Z_n$ corresponds to the Bernoulli's
 random walk), then functional (\ref{431}) can be represented at the form
\be\label{451} \tilde \phi^{s,t}_n={1\over \sqrt
n}\#\figbrakes{k\in[sn,tn): Z_n(k)=0}. \ee The functional
(\ref{451}) is widely used in a literature as the difference
analogue of the local time at the point zero for \emph{lattice}
random walks. Proposition \ref{c42} shows that the functional
(\ref{431}) is a natural difference analogue of the  local time
both for  random walks and, more generally, for difference
approximations of diffusion processes without any restrictions on
the  distribution of the sequence $\{\xi_k\}$.
\end{rem}

\section{Invariance principle for additive functionals of Markov chains}

In previous two chapters we have considered  more or less
particular examples illustrating possible ways to provide  the
main condition of Theorem \ref{t31} (condition 2). In this chapter
we introduce general sufficient condition of weak convergence of
additive functionals, constructed on the sequence of Markov
chains, that is formulated in terms of the transition
probabilities of these chains and the functions $F_n$ involved in
representation (\ref{31}). This condition is obtained as an
application of Theorem \ref{t31}, and the main assumption here is
that the local limit theorem (condition 4 of Theorem \ref{t41}
below) takes place in an appropriate form. For \emph{recurrent}
Markov chains this condition, together with a natural condition of
weak convergence of "symbols"\phantom{} of additive functionals
(exact formulation is given below), is sufficient for convergence
of characteristics, and  the estimates here are similar to
(\ref{42}) -- (\ref{43}) (see Theorem \ref{t62} below). For
\emph{transient} chains these estimates are not powerful enough,
since in this case the estimate (\ref{4_estimate}) does not
provide that $\Xi_n$ is negligible. One possible way to overcome
this difficultly is to apply a more strong version of local limit
theorem, for instance, to claim explicitly the rate of convergence
$\eps_k\to 0$ in (\ref{41}). Such an approach restricts the range
of possible applications, therefore we introduce another one, that
is concerned  with a uniform condition on the modulus of
continuity of processes $X_n$ (condition 5 of Theorem \ref{t41})
 and a "dimensional"\phantom{}  condition on the symbols of
 functionals (condition 6), adjusted one with another with an appropriate way
 (condition 7).

 We assume that a $\sigma$-finite
measures $\nu,\nu_n$ on $\Xf$ are given such that
$$
P(X(t)\in dy|X(s)=x)=p_{t-s}(x,y)\nu(dy),\quad 0\leq
s<t,x,y\in\Xf,
$$
$$
P\brakes{X_n\brakes{{i+k\over n}}\in dy|X_n\brakes{{i\over
n}}=x}=p_{n,k}(x,y)\nu_n(dy),\quad i\in \ZZ_+, k\in \NN,x,y\in\Xf.
$$
The measurable functions $p_t, p_{n,k}$ are interpreted as the
transition probability densities for $X,X_n$ w.r.t.  measures
$\nu,\nu_n$.

We assume the $W$-functional $\phi=\phi(X)$ with the
characteristics $f$ to be given. It is known (see \cite{dynkin},
Chapter 6) that
$$\phi^{s,t}=L_2-\lim_{\eps\to 0+}\int_s^t {1\over
\eps}f^{0,\eps}(X(r))\,dr,$$ and therefore
$$f^{s,t}(x)=\lim_{\eps\to 0+}\int_s^t\int_{\Xf}p_r(x,y) {1\over
\eps}f^{0,\eps}(y)\nu(dy)\,dr.
$$
We assume that, as $\eps\to 0+$, the measures  ${1\over
\eps}f^{0,\eps}d\nu$ converge weakly (i.e., on every bounded
continuous function)  to a  finite measure $\mu$,  the
characteristics $f$ can be represented in the form
\be\label{5_charact}
f^{s,t}(x)=\int_0^{t-s}\int_{\Xf}p_r(x,y)\mu(dy)\,dr,\quad \hbox{
and }\quad
\int_0^T\left[\sup_{x\in\Xf}\int_{\Xf}p_r(x,y)\mu(dy)\right]\,dr<+\infty,\quad
T\in \ax. \ee

We also consider the sequence of the functionals
$\phi_n=\phi_n(X_n)$ of the type (\ref{31}) with $L=1$  and
$F_{n}={1\over n}g_n$ (the case $L>1$ can be considered similarly
and we omit it in order to shorten notation). The characteristics
of $\phi_n$ has the form
$$
f_n^{s,t}(x)={1\over n}\sum_{s\leq {k\over
n}<t}\int_{\Xf}p_{n,k}(x,y)\mu_n(dy),\quad 0\leq s<t, x\in \Xf,
$$
where $\mu_n(dy)\equiv g_n(y)\nu_n(dy)$ are the
"symbols"\phantom{} of the functionals $\phi_n$.

\begin{thm}\label{t41} Assume the following conditions to hold true.
\begin{enumerate}\item Trajectories of the processes $X_n$ are
continuous,
 and the sequence $\{X_n\}$  possesses Markov approximation of $X$.

\item ${1\over n}\sup_x g_n(x)\to 0,n\to+\infty.$

\item For arbitrary $t_0>0$, the  function $(t,x,y)\mapsto
p_t(x,y)$ is uniformly continuous on $[t_0,+\infty)\times\Xf^2$,
and for arbitrary $y\in \Xf$
$$
\sup_{x\not\in B(y,R)}p_t(x,y)\to 0, \quad R\to+\infty
$$
(here and below $B(x,R)\equiv \{x\in\Xf|\rho(x,y)<R\}$).
Furthermore, there exist constants $\gamma>0, C_\gamma>0$ such
that
$$\sup_{x,y\in\Xf}p_{t}(x,y)\leq
C_\gamma t^{-\gamma}, \quad t>0.
$$

\item There exist sequences $\{\alpha_n\},\{\beta_n\}\subset
\Re^+$ tending to zero, such that
$$\sup_{x,y\in\Xf}|p_{n,k}(x,y)-p_{k\over n}(x,y)|\leq
(\alpha_n+\beta_k)\brakes{n\over k}^\gamma, \quad n,k\in \NN.
$$

\item There exist constants $\delta>0,C_\delta>0$ such that, for
arbitrary $T>0$,
$$
\sup_{x\in\Xf,n\in\NN} E\left(\left[\sup_{t,s\in[0,T],|t-s|\geq
{1\over n}}\,{\rho(X_n(t),X_n(s))\over
|t-s|^\delta}\right]^{C_\delta} \Bigl |X(0)=x\right)<+\infty.
$$

\item  Measures $\mu_n$ are finite and converge weakly to measure
$\mu$. There exist constants $\theta>0, C_\theta,c_\theta>0$ such
that
$$\mu_n(B(x,R))\leq C_\theta
R^\theta,\quad x\in\Xf,n\in \NN, R>c_\theta n^{-\delta}
$$
(note that the latter condition provides that  $\mu(B(x,R))\leq
C_\theta R^\theta, x\in\Xf,R>0$).

 \item The constants $\gamma,\delta,\theta,C_\delta$ satisfy the
 relations
$$
\delta\theta+1>\gamma, \quad C_\delta>2\theta+2.
$$
\end{enumerate}

Then $(X_n,\psi_n(X_n))\Rightarrow(X,\phi(X))$ in a sense of
convergence in distribution in $C(\ax,\Xf)\times C(\TT,\Re^+)$
($\psi_n$ are the random broken lines corresponding to the
functionals $\phi_n$).
\end{thm}

{\it Proof.} In order to prove the Theorem, it is sufficient to
show that, for every $T\in\ax$, \be\label{53}
f_n^{s,t}(x)\mathop{\rightrightarrows}\limits_{s\leq t\leq
T,x\in\Xf}f^{s,t}(x),\quad n\to+\infty. \ee Indeed, the sequence
$\{X_n\}$  provides Markov approximation for $X$ (condition 1),
and condition 1 of Theorem \ref{t31} is provides by condition 2 of
Theorem \ref{t41}. Having (\ref{53}) proved, we provide condition
2 of Theorem \ref{t31}. Condition 3 of this theorem is provided by
(\ref{5_charact}) and uniform continuity of the density $p$. At
last, condition 5 of Theorem \ref{t41} provides weak convergence
of $X_n$ to $X$ in $C(\ax,\Xf)$, that allows one to apply Theorem
\ref{t31} and Remark \ref{r32}.

Before proving (\ref{53}), let us  make some auxiliary estimates.
Denote
$$
H_{\delta,n}^{s,t}(X_n)=\sup_{v,w\in[s,t],|v-w|\geq {1\over
n}}\,{\rho(X_n(v),X_n(w))\over |v-w|^\delta},$$
$$
D_{n,A}^{s,t}=\figbrakes{X_n(r)\in
B(X_n(s),A(r-s)^\delta),r\in\sqrbrakes{s+{1\over n},t}},
$$
note that $\{H_{\delta,n}^{s,t}(X_n)<A\}\subset D_{n,A}^{s,t}$.
Also denote
$\alpha=\max_n\alpha_n,\beta=\max_k\beta_k,\delta_n={\sup_x|g_n(x)|\over
n}$, $B_1=\max_n\delta_n,
B_2(T)={C_\theta(C_\gamma+\alpha+\beta)\over
1+\delta\theta-\gamma}T^{1+\delta\theta-\gamma}.$ For arbitrary
$A>c_\theta,T\in\ax$, consider the functionals
 $\phi_{n,A}^{s,t}=\phi_n^{s,t}\1_{H_{\delta,n}^{s,t}(X_n)<A},s\leq t\leq T$.

 \begin{lem}\label{l51} 1. $E\sqrbrakes{\phi_{n,A}^{s,t}|X_n(s)=x}\leq B_1+B_2(T)A^\theta.$

 2. $E\sqrbrakes{\brakes{\phi_{n,A}^{s,t}}^2|X_n(s)=x}\leq {3B_1}(B_1+B_2(T)A^\theta)+2(B_1+B_2(T)A^\theta)^2$.

3.  Let $p\in\brakes{1,{2C_\delta-2\over C_\delta+2\theta}}$
(recall that $1<{2C_\delta-2\over C_\delta+2\theta}$ due to
condition 7 of the Theorem). Then
$$\sup_{x\in \Xf,n\in\NN,s\leq t\leq T}
E\sqrbrakes{\brakes{\phi_n^{s,t}}^p|X_n(s)=x}<+\infty. $$
 \end{lem}

{\it Proof.} Using condition 4 of the Theorem and then condition
6,    we obtain, for $t,s\in {1\over n}\ZZ_+$, the  estimate
$$ E\sqrbrakes{\phi_{n,A}^{s,t}|X_n(s)=x}\leq
E\sqrbrakes{\phi_{n}^{s,t}\1_{D_{n,A}^{s,t}}|X_n(s)=x}={g_n(x)\over
n}+
$$
$$
+{1\over n}\sum_{k=1}^{n(t-s)-1}\int_{B\brakes{x,A\brakes{k\over
n}^\delta}}p_{n,k}(x,y)\mu_n(dy)\leq
$$
$$
\leq \delta_n+{C_\gamma+\alpha+\beta\over n}\sum_{k=1}^{n(t-s)-1}
\brakes{n\over k}^\gamma \mu_n\brakes{B\brakes{x,A\brakes{k\over
n}^\delta}}\leq
 \delta_n+{C_\gamma+\alpha+\beta\over
n}\sum_{k=1}^{n(t-s)-1} C_\theta A^{\theta}\brakes{n\over
k}^{\gamma-\delta\theta}, $$ that immediately proves the first
statement of the Lemma. The second statement can be obtained from
the first one via the estimate similar to (\ref{391}) with the use
of the inequality
$$
\1_{H_{\delta,n}^{s,t}(X_n)<A}\leq
\1_{H_{\delta,n}^{s,r}(X_n)<A}\1_{H_{\delta,n}^{r,t}(X_n)<A},
$$
that holds true for arbitrary $r\in (s,t)$.

 Applying statement 2 and H\"older inequality we obtain
$$
E\sqrbrakes{\brakes{\phi_n^{s,t}}^p|X_n(s)=x}=\sum_{N=1}^\infty
E\sqrbrakes{\brakes{\phi_n^{s,t}}^p\1_{H_{\delta,n}^T(X_n)\in[N-1,N)}|X_n(s)=x}\leq
$$
$$
\leq \sum_{N=1}^\infty
E\sqrbrakes{(\phi_n^{s,t})^2\1_{H_{\delta,n}^T(X_n)<N}|X_n(s)=x}^{p\over
2}\sqrbrakes{P(H_{\delta,n}^T(X_n)\geq N-1)}^{2-p\over 2}\leq
$$
$$
\leq \sum_{N=1}^\infty
\sqrbrakes{B_3(T)+B_4(T)N^{2\theta}}^{p\over 2}\cdot
B_5(T)\sqrbrakes{(N-1)\vee 1}^{-{2-p\over 2} C_\delta},
$$
here and below $B_i(T), i=3,4,\dots$ denotes a constant, that can
be expressed explicitly through $T$ and the constants introduced
in the formulation of the Theorem, but an explicit expression is
not needed in our consideration.  Since ${\theta p}-{2-p\over
2}C_\delta <-1$ by the choice of $p$, this proves the statement 3.
The lemma is proved.

Let us proceed with the proof of (\ref{53}). Choose non-increasing
Lipschitz function $\Psi:\ax\to[0,1]$ such that
$\Psi([0,1])=\{1\},\Psi([2,+\infty))=\{0\}$, and set
$$
\Psi_r(x,y)=\Psi(r^{-1}\cdot \rho(x,y)),\quad r>0,x,y\in \Xf,\quad
\Psi_0\equiv 1.
$$
Note that, for arbitrary $r_0>0$, the function $(r,x,y)\mapsto
\Psi_r(x,y)$ is uniformly continuous on $[r_0,+\infty)\times
\Xf^2$.

For fixed $s\leq t\leq T, A\in\ax$ we decompose  $\phi_n^{s,t}$ as
$\phi_n^{s,t}=\eta_{n,A}^{s,t}+\zeta_{n,A}^{s,t},$ where
$$
\eta_{n,A}^{s,t}={1\over n}\sum_{s\leq {k\over
n}<t}g_n\brakes{X_n\brakes{k\over n}}\Psi_{A\brakes{{k\over n}-s
}^\delta}\left(X_n(s),X_n\brakes{k\over n}\right).
$$
We have that, on the set $D_{n,A}^{s,t}$, for $k$ such that $s\leq
{k\over n}<t$,
$$
\rho\brakes{X_n(0),X_n\brakes{k\over n}}\leq A\brakes{{k\over
n}-s}^\delta \Rightarrow \Psi_{A\brakes{{k\over n}-s
}^\delta}\left(X_n(s),X_n\brakes{k\over n}\right)=1,
$$
hence $\{\phi_n^{s,t}=\eta_{n,A}^{s,t}\}\supset D_{n,A}^{s,t}$ and
\be\label{54}\{\zeta_{n,A}^{s,t}\not=0\}\subset \Omega\backslash
D_{n,A}^{s,t}\subset\{H_{\delta,n}^{s,t}\geq A\}.\ee Let $p$ be
the same as in statement 3 of Lemma \ref{l51}. Then it follows
from (\ref{54}) and inequality $0\leq \zeta_{n,A}^{s,t}\leq
\phi_n^{s,t}$ that  \be\label{55}
E\sqrbrakes{\zeta_{n,A}^{s,t}|X_n(s)=x}\leq
E\sqrbrakes{(\phi_n^{s,t})^p|X_n(s)=x}^{1\over
p}\sqrbrakes{P(H_{\delta,n}^{s,t}\geq A|X_n(s)=x)}^{p-1\over
p}\leq B_6(T) A^{-\delta{p-1\over p}}. \ee
 Similarly, one can write
 $\phi^{s,t}=\eta_{A}^{s,t}+\zeta^{s,t}_A$, where
$ \eta_A^{s,t}=\int_s^t\Psi_{A(r-s)^\delta}(X(s),X(r))d\phi^{s,r},
$ \be\label{56} E\sqrbrakes{\zeta_{A}^{s,t}|X(s)=x}\leq  B_6(T)
A^{-\delta{p-1\over p}}. \ee

We have
$$
\abs{E\sqrbrakes{\eta_{n,A}^{s,t}|X_n(s)=x}-E\sqrbrakes{\eta_{A}^{s,t}|X(s)=x}}=
$$
$$=
\left|{g_n(x)\over n}+{1\over
n}\sum_{k=1}^{]n(t-s)[-1}\int_{\Xf}p_{k,n}(x,y)\Psi_{A\brakes{k\over
n}^\delta}(x,y)\mu_n(dy)-\int_0^{t-s}\int_{\Xf}p_r(x,y)\Psi_{Ar^\delta}(x,y)\mu(dy)\,dr\right|\leq
$$
$$
\leq
\delta_n+\Delta_n^1(x,A,s,t)+\Delta_n^2(x,A,s,t)+\Delta_n^3(x,A,s,t),
$$
where $]z[\equiv \min\{N\in\ZZ, N\geq z\}$,
$$
\Delta_n^1(x,A,s,t)=\left|{1\over
n}\sum_{k=1}^{]n(t-s)[-1}\int_{\Xf}[p_{k,n}(x,y)-p_{k\over
n}(x,y)]\Psi_{A\brakes{k\over n}^\delta}(x,y)\mu_n(dy)\right|,
$$
$$
\Delta_n^2(x,A,s,t)=\left|{1\over
n}\sum_{k=1}^{]n(t-s)[-1}\int_{\Xf}p_{k\over
n}(x,y)\Psi_{A\brakes{k\over
n}^\delta}(x,y)\mu_n(dy)-\int_0^{t-s}\int_{\Xf}p_r(x,y)\Psi_{Ar^\delta}(x,y)\mu_n(dy)\,dr\right|,
$$
$$
\Delta_n^3(x,A,s,t)=\left|\int_0^{t-s}\int_{\Xf}p_r(x,y)\Psi_{Ar^\delta}(x,y)[\mu_n(dy)-\mu(dy)]\,dr\right|.
$$
Denote $\Delta_n^i(A,T)= \sup_{x\in\Xf,s\leq t\leq
T}\Delta_n^i(x,A,s,t), i=1,2,3.$  Since $\Psi_r(x,y)\in [0,1]$ and
$\{\Psi_r(x,y)\not=0\}\subset \{y\in B(x,2r)\}$,
$$
\Delta_n^1(A,T)\leq {1\over
n}\sum_{k=1}^{]nT[-1}[\alpha_n+\beta_k]\brakes{n\over
k}^\gamma\mu_n\brakes{B\brakes{x,2A\brakes{k\over n}^\delta}}\leq
$$
\be\label{57}\leq  C_\theta(2A)^\theta\cdot {1\over
n}\sum_{k=1}^{]nT[-1}[\alpha_n+\beta_k]\brakes{k\over
n}^{\delta\theta-\gamma}\to 0,\quad n\to+\infty \ee by Toeplitz
theorem.

The function $(r,x,y)\mapsto p_r(x,y)\Psi_r(x,y)$ is uniformly
continuous over $[r_0,+\infty)\times X^2$ for any $r_0>0$,
therefore an estimate analogous to (\ref{57}) provides that
$$
\sup_{x\in\Xf,s\leq t\leq T}\left|{1\over
n}\sum_{k=[r_0n]+1}^{]n(t-s)[-1}\int_{\Xf}p_{k\over
n}(x,y)\Psi_{A\brakes{k\over
n}^\delta}(x,y)\mu_n(dy)-\int_{r_0}^{t-s}
\int_{\Xf}p_r(x,y)\Psi_{Ar^\delta}(x,y)\mu_n(dy)\,dr\right|\to
0
$$
(note that $\max_n\mu_n(\Xf)<+\infty$ since $\mu_n$ weakly
converge to $\mu$). The same arguments provide that
$$\mathop{\lim\sup}\limits_{n\to+\infty}\Delta_n^2(A,T)\leq
$$
$$
 \leq\mathop{\lim\sup}\limits_{n\to+\infty}\left[{1\over n}\sum_{k=1}^{[r_0n]}C_\gamma\brakes{n\over
k}^\gamma C_\theta\brakes{2A\brakes{k\over
n}^\delta}^{\theta}+\int_0^{r_0}C_\gamma
r^{-\gamma}C_\theta\brakes{2Ar^\delta}^\theta\,dr\right]=B_7(A,T)
(r_0)^{\delta\theta-\gamma+1}.
$$
Since $r_0>0$ is arbitrary, this implies that \be\label{58}
\Delta_n^2(A,T)\to 0,\quad n\to +\infty.\ee At last, the weak
convergence of $\mu_n$ to $\mu$ and the first part of condition 3
provide that, for every $t$,
$$
I_{n}(A,t)\equiv
\sup_{x\in\Xf}\left|\int_{\Xf}p_t(x,y)\Psi_{Ar^\delta}(x,y)[\mu_n(dy)-\mu(dy)]\right|\to
0,\quad n\to +\infty.
$$
Since $I_{n}(A,t)\leq C_\gamma t^{-\gamma}\cdot
C_\theta(2At^\delta)^\theta$, the  Lebesgue theorem of dominated
convergence provides  that \be\label{59} \Delta_n^3(A,T)\to
0,\quad n\to +\infty.\ee It follows from the estimates (\ref{55})
-- (\ref{59})  that
$$
\mathop{\lim\sup}\limits_{n\to+\infty}\sup_{x\in\Xf,s\leq t\leq
T}\abs{f_n^{s,t}(x)-f^{s,t}(x)}\leq 2B_6(T)A^{-\delta{p-1\over
p}},\quad A>c_\theta.
$$
Taking $A\to+\infty$ we obtain (\ref{53}), that completes  the
proof. The theorem is proved.

In order to make our exposition complete, let us formulate a
version of Theorem \ref{t41} for the recurrent case.

\begin{thm}\label{t62} Let conditions 1 -- 4 of Theorem
\ref{t41} hold true and $\gamma<1$. Also let $\mu_n$ converge
weakly to $\mu$, and $X_n$ converge to $X$ by distribution in
$C(\ax,\Xf)$.

Then $(X_n,\psi_n(X_n))\Rightarrow(X,\phi(X))$ in a sense of
convergence in distribution in $C(\ax,\Xf)\times C(\TT,\Re^+)$.
\end{thm}

The proof, with slight changes, repeats the proof of Theorem
\ref{t41}, and is omitted. Note that, under conditions of Theorem
\ref{t62}, the convergence of finite-dimensional distributions of
$\phi_n$ can be provided with the use of the technique, mentioned
in the Introduction, that was proposed by I.I.Gikhman and is based
on studying of limit behavior of difference equations for
characteristic functions of $\phi_n^{s,t}$ (see for instance the
proof of Theorem 3 \cite{portenko_68}). In the transient case,
treated in Theorem \ref{t41}, this technique can not be applied
since the uniform estimates, analogous to (\ref{42}) --
(\ref{43}), are not available in this case.

At last, let us give an example of application of Theorem
\ref{t41}. To shorten exposition we omit the proofs of some
technical details.

\begin{ex} Let $\Xf=\Re^d,d\geq 2$ and $X_n,X$ be as in
Example 1. Let $K\subset \Re^d$ be a compact set, for which the
surface measure $\lambda_K$ is well defined by equality
$$
\lambda_K(\cdot)\equiv w-\lim_{\eps\to 0+}{\lambda^d(\cdot\cap
K_\eps)\over \lambda^d(K_\eps)},
$$
where $w-\lim$ means the limit in the sense of weak convergence of
measures, $\lambda^d$ is Lebesgue measure on $\Re^d$,
$K_\eps\equiv \{x|\mathrm{dist}(x,K)\leq \eps\}$. Assume that the
condition \be\label{61} \lambda^d(K_\eps)\geq \mathrm{const}\cdot
\eps^{\beta},\quad \eps>0\ee holds with some $\beta<2$. In
particular, the set $K$ can be smooth (or, more generally,
Lipschitz) surface of codimension 1 or fractal with its
Haussdorf-Besikovich dimension greater then $d-2$.

It not hard to verify that $\mu\equiv\lambda_K$ is $W$-measure
(see \cite{dynkin}, Chapter 8.1 for the terminology), and
therefore corresponds to some $W$-functional $\phi$ of the Wiener
process $X$. This functional is naturally interpreted as the local
time of Wiener process at the set $K$, and can be written as
$\phi^{s,t}=\int_s^t \lambda_K(X_r)\,dr.$

We consider the functionals $\phi_n(X_n)$ of the form
$$
\phi_s^{t,n}={1\over n\lambda^d\brakes{K_{1\over {\sqrt
n}}}}\sum_{k\in[sn,tn)}\1_{\{X^n({k\over n})\in K_{1\over {\sqrt
n}}\}},
$$
and apply Theorem \ref{t41} in order to prove convergence of the
distributions in $C(\ax,\Re^d)\times C(\TT,\Re^+)$
\be\label{60}(X_n,\psi_n(X_n))\Rightarrow(X,\phi(X))\ee ($\psi_n$
are the broken lines corresponding to $\phi_n$).

 Condition 1 holds true due to Example 1, condition 2 is provided
 by condition (\ref{61}) (by this condition,
 $\sup_x g_n(x)\leq \mathrm{const}\cdot n^{\beta\over 2}$).
Condition 3 holds with $p_t(x,y)=(2\pi t)^{-{d\over
2}}\exp\{-{1\over 2}\|y-x\|^2_{\Re^d}\}$ and $\gamma={d\over 2}$.
Condition (\ref{61}) implies condition 6  with $\theta=d-\beta.$

We assume that the random walk $S_n$ is either aperiodic on some
lattice $h\ZZ^d$ or is strongly non-lattice (i.e., $S_{n_0}$ has
bounded distribution density for some $n_0$). Under this
assumption, condition 4 holds with $\alpha_n\equiv 0$,
$\nu=\lambda^d$ and $\nu_n$ equal to counting measures on ${h\over
\sqrt n}\ZZ^d$ in lattice case or $\lambda^d$ in strongly
non-lattice case.

It remains to provide conditions 5, 7. We have ${\gamma-1\over
\theta}={d-2\over 2(d-\beta)}<{1\over 2}$. Choose some $\delta\in
\left({\gamma-1\over \theta},{1\over 2}\right)$ and consider
$\alpha>0$ such that ${{\alpha\over 2}-1\over \alpha}>\delta$ and
$\alpha>2\theta+2$. Suppose that  \be\label{62}
E\|\xi_k\|_{\Re^d}^\alpha<+\infty. \ee Then applying Burkholder
inequality we obtain that
\be\label{63}E\|X_n(t)-X_n(s)\|^\alpha_{\Re^d}\leq
\mathrm{const}\cdot |t-s|^{\alpha\over 2},\quad |t-s|\geq{1\over
\sqrt n}, x\in\Re^d.\ee  Repeating the standard proof of the
Kolmogorov's theorem on existence of continuous modification (see,
for instance \cite{skor_tsp}, p. 44,45), one can deduce from
(\ref{63}) that, for $\varsigma<\alpha, \vartheta<{{\alpha\over
2}-1\over \alpha}$,
$$
\sup_n E\left[\sup_{t,s\in[0,T],|t-s|\geq {1\over
n}}\,{\|X_n(t)-X_n(s)\|_{\Re^d}\over
|t-s|^\vartheta}\right]^{\varsigma} <+\infty.
$$
Finally, choosing $\vartheta=\delta,\varsigma>2\theta+2$ we obtain
that conditions 5,7 hold with $C_\theta=\varsigma$. Applying
Theorem \ref{t41}, we obtain  weak convergence (\ref{60}) under
additional moment condition (\ref{62}). One can remove this
condition using the "cutting" procedure, described in the proof of
the Proposition \ref{c42}.

Let us remark that for the lattice random walks the result,
exposed in Example 4, was obtained in \cite{bass_khoshn92} by a
technique, essentially different from the one proposed here.
Convergence (\ref{60}) in continuous case, as far as it is known
to authors, is a new result.

\end{ex}

\end{document}